\documentclass[a4paper,11pt]{article}
\usepackage{amsmath}
\usepackage{amsthm}
\usepackage{amssymb}
\usepackage{amsfonts}
\usepackage{enumerate}

\normalsize

\newtheoremstyle{theoremstyle}
  {10pt}      %  Space above
  {5pt}       %  Space below
  {\itshape}  %  Body font
  {}          %  Indent amount (empty = no indent, \parindent = para indent)
  {\bfseries} %  Thm head font
  {:}         %  Punctuation after thm head
  {.5em}      %  Space after thm head: " " = normal interword space;
  {}          %  Thm head spec (can be left empty, meaning `normal')

\newtheoremstyle{examplestyle}
  {10pt}      %  Space above
  {5pt}       %  Space below
  {}          %  Body font
  {}          %  Indent amount (empty = no indent, \parindent = para indent)
  {\bfseries} %  Thm head font
  {:}         %  Punctuation after thm head
  {.5em}      %  Space after thm head: " " = normal interword space;
  {}          %  Thm head spec (can be left empty, meaning `normal')

\theoremstyle{theoremstyle}
\newtheorem{theorem}{Theorem}[section]
\newtheorem*{theorem*}{Theorem}
\newtheorem{lemma}[theorem]{Lemma}
\newtheorem{proposition}[theorem]{Proposition}
\newtheorem*{proposition*}{Proposition}
\newtheorem{corollary}[theorem]{Corollary}
\newtheorem*{corollary*}{Corollary}

\theoremstyle{examplestyle}
\newtheorem{example}[theorem]{Example}
\newtheorem{definition}[theorem]{Definition}
\newtheorem{definition*}{Definition}

\newtheorem{remark*}{Remark}

\newcommand{\comment}[1]{}

\newcommand{\rays}{{\Delta(1)}}

\newcommand{\sh}[1]{\mathcal{#1}}

\newcommand{\spec}{\operatorname{spec}}

\newcommand{\orb}{\operatorname{orb}}

\newcommand{\codim}{\operatorname{codim}}

\newcommand{\rk}{\operatorname{rk}}

\newcommand{\Hom}{\operatorname{Hom}}
\newcommand{\Ext}{\operatorname{Ext}}
\newcommand{\coker}{\operatorname{coker}}

\newcommand{\ksm}{{k[\sigma_M]}}

\newcommand{\supp}{{\operatorname{supp}}}

\newcommand{\depth}{{\operatorname{depth}}}
\newcommand{\stern}{{\operatorname{star}}}

\newcommand{\Z}{\mathbb{Z}}
\newcommand{\R}{\mathbb{R}}

\title{On the Local Cohomology of Reflexive Modules of Rank One over Normal
Semigroup Rings}

\author{Markus Perling\footnote{Institut f\"ur Mathematik, Universit\"at
Paderborn, 33098 Paderborn, Germany, {\tt perling@math.upb.de}}}

\date{March 24, 2005}

\begin{document}

\maketitle

\begin{abstract}
In this work we describe the local cohomology of reflexive modules of rank one over
normal semigroup rings with respect to monomial ideals. Using our description we show
that the problem of classifying maximal Cohen-Macaulay modules of rank one can be
rephrased in terms of finding integral solutions to certain sets of linear
inequalities.
\end{abstract}

\section{Introduction}

The motivation for this paper stems from my earlier and ongoing work on
equivariant sheaves over toric varieties (see \cite{perling1}, \cite{perlingdiss}
and \cite{perling4}). The main theme of this theory is the interplay between the
combinatorics of toric geometry and non-combinatorial aspects from linear algebra.
In a sense, this theory extends the combinatorial theory of toric varieties to a
semi-combinatorial theory {\em over} toric varieties. It turns
out that an important building block which we should understand are the reflexive
sheaves of rank one.  This paper has been written in order to clarify at least
a few aspects of these sheaves.

 In the case of affine toric varieties, these sheaves correspond to
reflexive modules of rank one over a normal semigroup rings $\ksm$, where $k$ is an
algebraically closed field and $\sigma_M$ a normal subsemigroup of some lattice
$M \cong \Z^d$. The main results of this work are related to the following problems:
\begin{enumerate}[(i)]
\setlength{\itemsep}{-5pt}
\item the computation of local cohomology modules of reflexive modules of rank one over
$\ksm$ with respect to monomial ideals,
\item the classification of maximal Cohen-Macaulay (MCM) modules of rank one over
$\ksm$.
\end{enumerate}

 Recall that a normal semigroup is of the form
\begin{equation*}
\sigma_M = \{m \in M \mid \langle m, n(\rho) \rangle \geq 0 \text{ for all }\, \rho
\in \sigma(1)\},
\end{equation*}
where the $n(\rho)$ are linear forms in $N = M\check{\ }$, the module dual to $M$, and
the brackets $\langle\,  , \, \rangle : M \times N \rightarrow \Z$ denote the canonical
pairing. The $n(\rho)$ over $\R_{\geq 0}$ span a strictly convex polyhedral
cone $\sigma$ in the vector space $N_\R = N \otimes_Z \R$ (see also section
\ref{prelim}). Denote $T = \Hom(M, k^*)$ the algebraic torus acting on the affine toric
variety $U_\sigma$. The $T$-invariant divisor class group is isomorphic to the free
group $\Z^{\sigma(1)}$. By a well-known correspondence, there is a one-to-one
correspondence between $T$-invariant divisors $D \in \Z^{\sigma(1)}$ and $M$-graded
reflexive module
of rank one over $\ksm$, denoted $R^D$. Forgetting the grading on $R^D$ in a natural
way induces a one-to-one correspondence between isomorphism classes of reflexive
modules of rank one and rational equivalence classes of Weil-divisors.
Let $D = \sum_{\rho \in \sigma(1)} n_\rho D_\rho$, then, as an $M$-graded module,
$R^D$ can explicitly be represented as the $k$-linear span of a $\sigma_M$-lattice
submodule $M^D$ of $M$, which is determined by inequalities
\begin{equation*}
M^D = \{m \in M \mid \langle m, n(\rho) \rangle \geq -n_\rho \text{ for all }\, \rho
\in \sigma(1)\}.
\end{equation*}

\paragraph{Local cohomology.} Using this representation of $R^D$, the main technical
result of this paper will be a characterization in terms of simplicial cohomology of
the local cohomology of $R^D$ with respect to a monomial ideal. For any ring $S$, an
ideal $B$ of $S$ and any $S$-module $F$, there is the functor
\begin{equation*}
\Gamma_B F := \underset{\underset{n}{\rightarrow}}{\lim} \Hom(S / B^n, F).
\end{equation*}
The local cohomology modules $H^i_B F$ are defined as the right derived functors of
$\Gamma_B$ and have the following characterization:
\begin{equation*}
H^i_B F = \underset{\underset{n}{\rightarrow}}{\lim} \Ext^i(S / B^n, F).
\end{equation*}
In general, an explicit description of the modules $H^i_B F$ is very difficult.
However, in our case $B$ is an ideal in $\ksm$ generated by monomials, and the
modules $H^i_B R^D$ are in a natural way $M$-graded and admit an explicit
combinatorial description. For this, we construct an $M$-graded resolution of $R^D$
(Proposition \ref{acycres}):
\begin{equation*}
\mathbf{D}^\cdot: 
0 \rightarrow R^D \rightarrow A_0 \longrightarrow A_1
\longrightarrow \cdots \longrightarrow A_{\vert\sigma(1)\vert} \longrightarrow 0,
\end{equation*}
where the $A_i$ are $\Gamma_B$-acyclic $\ksm$-modules. Applying $\Gamma_B$ to
$\mathbf{D}^\cdot$ then yields an isomorphism $H^i_B R^D \cong H^i(\Gamma_B
\mathbf{D}^\cdot)$.
To analyze this cohomology a bit deeper, we consider the {\em support} and
{\em cosupport} of $B$. The support of $B$ is defined to be set of faces $\tau$ of
$\sigma$ such that the corresponding orbit $\orb(\sigma)$ is contained in the variety
$V(B)$. To introduce the cosupport, we denote $\Sigma$ the simplex spanned by
$\sigma(1)$, i.e. the set underlying $\Sigma$ coincides with $\sigma(1)$, but we
choose an order
on this set, such that $\Sigma$ becomes an oriented combinatorial simplex. Then the
cosupport
$\Xi_B$ of $B$ is the set of those $\Pi \subset \Sigma$ such that the minimal face
$\tau$
of $\sigma$ where $\tau(1)$ contains $\Pi$, is not contained in the support of $B$.
Then $\Xi_B$ in a natural way can be considered as a simplicial subcomplex of $\Sigma$.
Now, for any $m \in M$, let $\Sigma_m := \{\rho \in \Sigma \mid \langle m, n(\rho)
\rangle < 0\}$. Then we have for the $m$-th graded component of
the complex $\Gamma_B \mathbf{D}^\cdot$ (Corollary \ref{cohomcorollary}):
\begin{equation*}
\big(H^i_B R^D\big)_m = \big(H^i(\Gamma_B \mathbf{D}^\cdot)\big)_m =
\tilde{H}^{i - 2}(\Xi_B \cap \Sigma_m; k);
\end{equation*}
i.e. we identify the graded components $\big(H^i_B R^D\big)_m$ with the $(i - 2)$-th
reduced cohomology of the simplicial complex $\Xi_B \cap \Sigma_m$. Note that this
kind of identification is not new but has been applied in the literature several times
to study the local cohomology of (not necessarily normal) semigroup rings (see, for
instance, \cite{GotoWatanabe}, \cite{TrungHoa}, \cite{Mustata2}). In fact, the
constructions in this work are a quite straightforward adaption of the methods of
Trung and Hoa
\cite{TrungHoa}. The new aspect here is that we apply this technique to the study of
more general $\ksm$-modules. We remark that local cohomology of general
$M$-graded modules has been studied before (see \cite{helmmiller}, \cite{helmmiller2}),
but here, as will be explained below, we will arrive
at a more explicit combinatorial picture for the case of the modules $R^D$.

\paragraph{Maximal Cohen-Macaulay modules of rank one.}

By a classical theorem of Hochster, the rings $\ksm$ are Cohen-Macaulay and it is a
natural problem to classify (maximal) Cohen-Macaulay modules over $\ksm$. Although
there exists a huge amount of literature concerning the classification of MCM modules
in various contexts, to my knowledge, there has not yet been done much work for
the case of rings $\ksm$. The only references I am aware of and which are explicitly
devoted to this topic, are \cite{BrunsGubeladze} and \cite{BrunsGubeladze2}.
One of the main results of \cite{BrunsGubeladze} (Corollary 5.2) --- from the
perspective of this paper, at least --- is that over a normal semigroup ring there
exist only finitely many isomorphism classes of MCM modules of rank one.

However, despite of this finiteness result, it seems that a complete classification
is very difficult, if not impossible, to achieve. The relevant combinatorics behind
such a classification is given by the hyperplane arrangement defined by the real
hyperplanes $H_\rho = \{m \in M \otimes_Z \R \mid \langle m, n(\rho) \rangle =
-n_\rho\}$ in $M \otimes_\Z \R$. To this hyperplane arrangement one
can associate its combinatorial type which is represented by the so-called
{\em matroid of flats}; this structure essentially captures the information on
intersections of the hyperplanes $H_\rho$.
It is a well-known fact that $R^D$ is MCM if and only if
$H^i_{\mathfrak{m}} R^D = 0$ for all $0 \leq i < d$, where $\mathfrak{m}$ is the
maximal $M$-graded ideal in $\ksm$. So, by the results above,
$R^D$ is MCM if and only if for every subset $\Pi \subset
\Sigma$ such that $\tilde{H}^{i - 2}(\Pi \cap \Xi_{\mathfrak{m}}, k) \neq 0$ for some
$i < d$, the
system of linear inequalities
\begin{align*}
\langle m, n(\rho) \rangle & < -n_\rho \quad \text{ for } \rho \in \Pi \\
\langle m, n(\rho) \rangle & \geq -n_\rho \quad \text{ for } \rho \in \Sigma
\setminus \Pi
\end{align*}
has no solution in $M$. So, we can phrase the problem of classifying MCM modules of
rank one over $\ksm$ as follows. Given the matroid of flats of a hyperplane
arrangement $\bigcup_{\rho \in \rays} H_\rho$, then, how many possibilities are there
to realize (up to translation) combinatorially equivalent hyperplane arrangements such
that a certain subset of the open cells of $M_\R \setminus \bigcup_{\rho \in \rays}
H_\rho$ (and part of their boundaries, respectively) does not intersect the lattice
$M$? A definitive general solution by now seems to be out of reach.

\paragraph{Overview of the paper.} In section \ref{prelim}, we introduce some general
facts and notation from toric geometry, which will be used throughout the paper. In
section \ref{sigmasimplex} we introduce the simplex over the set of rays of the fan
$\sigma$ and recall some elementary and basic facts on simplicial chain complexes. In
section \ref{acyclicsection} we construct a class of $\ksm$-modules which are
acyclic with respect to local cohomology functors of monomial ideals. We use these
modules to construct acyclic resolutions for the modules $R^D$ in section
\ref{resolution}. In that section we also characterize the graded components of local
cohomology modules in terms of reduced cohomology of certain cell complexes.
Section \ref{chambersection} presents some easy observations and examples concerning
the chambers of $M_\R \setminus \bigcup_{\rho \in \rays} H_\rho$ and the vanishing and
nonvanishing of local cohomology in certain degrees. Section \ref{mcmsection} is a
short insertion to give the clear statement about the conditions on $R^D$ to be MCM.
In section \ref{singsets} we collect some more facts about the depth of the modules
$R^D$.

\section{Toric preliminaries}
\label{prelim}

\paragraph{General notions.}
We introduce some notation from the theory of toric varieties.
For general overview on toric varieties we refer to \cite{Oda}, \cite{Fulton}.
We will always assume that $k$ is an algebraically closed field.
$U_\sigma$ will denote an affine toric variety over $k$ defined by a strongly convex
rational
polyhedral cone $\sigma$ contained in the real vector space $N_\mathbb{R}
\cong N \otimes_{\mathbb{Z}} \mathbb{R}$ over a lattice $N \cong \mathbb{Z}^d$.
We always assume that $\dim \sigma = \rk_\Z N$.
Let $M$ be the lattice dual to $N$ and let $\langle \, ,
\rangle : M \times N \rightarrow \mathbb{Z}$ be the canonical pairing. This pairing
extends in a natural way to the scalar extensions $M_\mathbb{R} := M
\otimes_\mathbb{Z} \mathbb{R}$ and $N_\mathbb{R}$. Elements of $M$ are denoted by
$m$, $m'$, etc. if written additively, and by $\chi(m)$, $\chi(m')$, etc. if written
multiplicatively, i.e. $\chi(m + m') = \chi(m)\chi(m')$. The lattice $M$ is identified
with the group of characters of the torus
$T = \operatorname{Hom}(M,k^*) \cong (k^*)^d$ acting on $U_\sigma$.
For any cone $\sigma$ we will use the following notation:

\begin{itemize}
\setlength{\itemsep}{-5pt}
\item faces of $\sigma$ are denoted by small Greek letters such as $\rho$,
$\tau$, etc., the natural order among faces is denoted by $\rho \prec \tau$;
\item for any set $F$ of faces of $\sigma$, we call the set
\begin{equation*}
\stern(F) := \{\eta \mid \tau \prec \eta \text{ for some } \tau \in F\}
\end{equation*}
the {\em star} of $F$; if $\stern(F) = F$, then $F$ is called {\em star closed};
\item $\sigma(1) := \{ \rho \prec \sigma \mid \dim \rho = 1\}$, the faces of dimension
1, called {\it rays}; below we will often find it more convenient to denote this set
$\Sigma$ (by abuse of notation);
\item $n(\rho)$ denotes the primitive lattice element spanning the ray $\rho$;
\item $D_\rho$ denotes the the irreducible $T$-invariant Weil divisor on $U_\sigma$
associated to $\rho \in \sigma(1)$;
\item $\check{\sigma} := \{m \in M_\mathbb{R} \mid \langle m, n(\rho) \rangle \geq 0
\text{ for all $\rho \in \sigma(1)$}\}$ is the cone {\it dual} to $\sigma$;
\item $\sigma^\bot = \{m \in M_\mathbb{R} \mid \langle m, n \rangle = 0 \text{ for
all } n \in \sigma \}$;
\item $\sigma_M := \check{\sigma} \cap M$ is the subsemigroup of $M$ associated 
to $\sigma$;
\item $\sigma_M^\bot := \sigma^\bot \cap M$ is the unique maximal subgroup of $M$
contained in $\sigma_M$;
\item the semigroup ring $\ksm \cong \bigoplus_{m \in \sigma_M} k \cdot \chi(m)$
is identified with the coordinate ring of $U_\sigma$, and the group ring $k[M]$
is identified with the coordinate ring of $T$;
\item let $\tau \prec \sigma$, then $\orb(\tau)$ denotes the orbit associated to $\tau$
in $U_\sigma$;
\item let $\tau \prec \sigma$, then $N_\tau$ is $N$ intersected with the subvector
space of $N$ spanned by $\tau$ over $\R$, $M_\tau$ is the dual module of $N_\tau$,
where there is a canonical identification $M_\tau = M / \tau^\bot_M$; moreover, we
there is the canonical splitting $M = \tau^\bot_M \times M_\tau$;
\item $T_\tau$ denotes the stabilizer subgroup of $T$ over $T_\tau$; $T^\tau := T /
T_\tau$, note that $T^\tau \cong \orb(\tau)$.
\end{itemize}
Recall that dualizing $\sigma$ via $\tau \mapsto \tau^\bot \cap \check{\sigma}$
induces an
order-reversing one-to-one correspondence between faces of $\sigma$ and faces of
$\check{\sigma}$.

\paragraph{Reflexive modules of rank one.}
There exists a short exact sequence of $\Z$-modules:
\begin{equation*}
0 \longrightarrow M \longrightarrow \Z^{\sigma(1)} \longrightarrow
A_{d - 1}(U_\sigma) \longrightarrow 0,
\end{equation*}
where $A_{d - 1}(U_\sigma)$ is the $(d - 1)$-st Chow group of $U_\sigma$, i.e. the
group of rational equivalence classes of Weil divisors on $U_\sigma$, and
$\Z^{\sigma(1)}$ the group which is freely generated over the $T$-invariant
irreducible Weil divisors of $U_\sigma$. It was observed by Reid \cite{Reid1} that
there is a one-to-one correspondence of classes $\alpha \in A_{n - 1}(U_\sigma)$
and isomorphism classes of reflexive sheaves $\sh{O}(D)$ of rank one over $U_\sigma$,
where $D$ is some representative for $\alpha$. Every reflexive module of rank one over
$\ksm$ is isomorphic to $\Gamma\big(U_\sigma, \sh{O}(D)\big)$ for some Weil divisor
$D$, and in fact, there is a one-to-one correspondence between rational equivalence
classes of Weil divisors and reflexive $\ksm$-modules of rank one.
The above short exact sequence implies that every class $\alpha$ can be
represented by a $T$-invariant Weil divisor $D = D_{\underline{n}}$, where
$\underline{n} = (n_\rho) \in \Z^{\sigma(1)}$ and $D = \sum_{\rho \in \sigma(1)}
n_\rho D_\rho$. By $T$-invariance, to $D$ there corresponds a reflexive $\ksm$-module
of rank one, denoted $R^D$, which has a natural $M$-graded structure together with
a natural $M$-graded embedding $R^D \hookrightarrow k[M]$. Namely,
denote $M^D_\rho := \{m \in M \vert \langle m, n(\rho) \rangle \geq -n_\rho\}$, which
corresponds to the shifted half space $\rho_M + m$, where $m \in M$ such that
$\langle m, n(\rho) \rangle = -n_\rho$. We define $M^D := \bigcap_{\rho \in
\sigma(1)} M^D_\rho$, then $R^D := k[M^D]$ is a well-defined $M$-graded reflexive
$\ksm$-submodule of $k[M]$.

\paragraph{Monomial ideals.}
Let $I \subset \sigma_M$ be a semigroup ideal, then $B := \bigoplus_{m \in I}
k \cdot \chi(m)$ is an ideal in $\ksm$. On the other hand, for every $M$-graded ideal
$B$ the set $I = \{m \mid B_m \neq 0\}$ forms a semigroup ideal in $\sigma_M$.
We call the class of ideals coming from semigroup ideals the {\em monomial ideals}.

\begin{definition}
Let $B$ be a monomial ideal, then its support is defined as
\begin{equation*}
\supp(B) := \{\tau \prec \sigma \mid I \cap \tau^\bot \cap \check{\sigma} = \emptyset\}.
\end{equation*}
\end{definition}

\noindent Note that $\supp(B)$ is star-closed and the variety $V(B)$ coincides with
$\bigcup_{\tau \in \supp(B)} \orb(\tau)$ $\subset U_\sigma$. Moreover, in the
particular case where $B$ is the unique maximal homogeneous ideal of $\ksm$, we have
$\supp(B) = \{\sigma\}$ and $V(B) = \orb(\sigma)$.

\section{The simplex spanned by $\sigma$}\label{sigmasimplex}

We denote $\Sigma$ the {\em simplex} spanned by $\sigma$, i.e. the set underlying
$\Sigma$ coincides with $\sigma(1)$, but we choose a total ordering on the elements
of $\Sigma$, i.e. $\Sigma = \{\rho_1 < \rho_2 < \cdots < \rho_n\}$. Any subset $\Pi
\subset \sigma(1)$ corresponds to a face of $\Sigma$ with orientation induced by the
ordering of $\Sigma$. In what follows we will find it convenient to identify
$\sigma(1)$ and $\Sigma$ as sets and to write $\Sigma$ instead of $\sigma(1)$.

For any $\Pi \subseteq \Sigma$ with $\vert \Pi \vert = r$, we consider the
augmented cochain complex:
\begin{equation*}
\mathbf{\Phi}_\Pi :
0 \longrightarrow \Z \overset{\Phi_{-1}}{\longrightarrow} \bigoplus_{\rho \in \Pi}
\Z \cdot \rho \overset{\Phi_0}{\longrightarrow} \bigoplus_{\substack{\Gamma \subset
\Pi\\ \vert \Gamma \vert = 2}} \Z \cdot \Gamma \overset{\Phi_1}{\longrightarrow}
\cdots \overset{\Phi_{r - 1}}{\longrightarrow} \Z \cdot \Pi \longrightarrow 0,
\end{equation*}
which is an exact sequence of $\Z$-modules, as $\Pi$ is contractible. To avoid to
mention the case $\Pi = \emptyset$ repeatedly as a special case in exactness
arguments, we adopt the convention that for $\Pi = \emptyset$ the corresponding
augmented
cochain complex is $0 \rightarrow \Z \rightarrow \Z \overset{\Phi_{-1}}{\rightarrow}
0$.
For any two subsets $\Pi \subset \Upsilon \subset \Sigma$, the canonical
projections
\begin{equation*}
\bigoplus_{\substack{\Gamma \subset \Upsilon\\ \vert \Gamma \vert = i}} \Z \cdot \Gamma
\twoheadrightarrow
\bigoplus_{\substack{\Gamma \subset \Pi\\ \vert \Gamma \vert = i}} \Z \cdot \Gamma
\end{equation*}
for every $i = 1, \dots, \vert \Upsilon \vert$, induce a surjective chain map
$\mathbf{\Phi}_\Upsilon \twoheadrightarrow \mathbf{\Phi}_\Pi$.

\begin{definition}
Let $\Pi \subset \Sigma$, then we denote
\begin{equation*}
\tau_\Pi := \min \{\tau \prec \sigma \mid \Pi \subset \tau(1)\},
\end{equation*}
i.e. $\tau_\Pi$ is the minimal face of $\sigma$ such that $\Pi \subset \tau(1)$.
If $\Pi = \emptyset$, then $\tau_\Pi$ is the zero cone.
\end{definition}

\begin{definition}
Let $F$ be any set of faces of $\sigma$, then we set
\begin{equation*}
\Xi_F := \{\Pi \subset \Sigma \mid \tau_\Pi \notin \stern(F)\},
\end{equation*}
which we consider as a simplical subcomplex of $\Sigma$.
If $F = \{\tau\}$ for some $\tau \prec \sigma$, then we write $\Xi_\tau$ instead of
$\Xi_{\{\tau\}}$, and we denote $\Xi_\sigma$ simply by $\Xi$.
If $B \subset \ksm$ is a monomial ideal and $F = \supp(B)$, then we also write
$\Xi_B$ instead of $\Xi_{\supp(B)}$. $\Xi_B$ is called the {\em cosupport} of $B$ in
$\Sigma$.
\end{definition}
Clearly, for every $\Upsilon \subset \Pi \subset \Sigma$, $\Pi \in \Xi_F$ implies
$\Upsilon \in \Xi_F$, and $\Xi_F$ is a simplicial subcomplex of $\Sigma$.

Now let $\Pi \subset \Sigma$ be any subset which we consider as subsimplex and let $F$
any set of faces of $\sigma$, then we set $\Pi \cap \Xi_F := \{\Upsilon \mid \Upsilon
\subset \Pi \text{ and } \Upsilon \in \Xi_F\}$, which is a subcomplex of both $\Pi$
and $\Xi_F$. Therefore there exists a subcomplex $\mathbf{\Phi}_{\Pi, \Xi_F}$ of
$\mathbf{\Phi}_\Pi$ which is build of the terms
\begin{equation*}
\mathbf{\Phi}_{\Pi, \Xi_F}^i := \bigoplus_{\substack{\Gamma \subset \Pi \cap \Xi_F\\
\vert\Pi\vert = i}} \Z \cdot \Gamma,
\end{equation*}
The relative reduced cohomology $\tilde{H}^\cdot(\Pi, \Pi \cap \Xi_F; \Z)$ then is
the homology of the quotient complex $\mathbf{\Phi}_\Pi / \mathbf{\Phi}_{\Pi, \Xi_F}$.
Note that because $\Pi$ is contractible, there are isomorphisms $\tilde{H}^i(\Pi, \Pi
\cap \Xi_F; \Z) \cong \tilde{H}^{i - 1}(\Pi \cap \Xi_F; \Z)$ for every $i \in \Z$.
Likewise, for our chosen field $k$, we obtain the relative reduced cohomology with
coefficients in $k$ as the homology of $\mathbf{\Phi}_\Pi \otimes_\Z k / \big(
\mathbf{\Phi}_{\Pi, \Xi_F} \otimes_\Z k\big)$ and $\tilde{H}^i(\Pi, \Pi \cap \Xi_F; k)
\cong \tilde{H}^{i - 1}(\Pi \cap \Xi_F; k)$ for every $i \in \Z$.

\section{$\Gamma_B$-Acyclic Modules}\label{acyclicsection}

\begin{definition}
Let $\Pi \subset \Sigma$ be a nonempty subset, then we set
\begin{equation*}
H^D_\Pi := M \setminus \bigcup_{\rho \in \Pi} M^D_\rho
\end{equation*}
and define
\begin{equation*}
k[H^D_\Pi] := \coker\big(\bigoplus_{\rho \in \Pi} k[M^D_\rho] \longrightarrow
k[M]\big).
\end{equation*}
In the special case where $\Pi = \{\rho\}$, we write $H^D_\rho$ and $k[H^D_\rho]$,
respectively.
\end{definition}
Clearly, we have $H^D_\Pi = \bigcap_{\rho \in \Pi} H^D_\rho$ and
$k[H^D_\Pi] = \bigoplus_{m \in H^D_\Pi} k \cdot \chi(m)$ a $\ksm$-module.
Let $\rho \in \Sigma$ and any $m \in H^D_\rho$, then for any subset $M' \subset
\rho^\bot_M$ the set $m + M'$ is contained in $H^D_\rho$. This is in particular true
when $M'$ is a subgroup of $\rho^\bot_M$. Thus the following lemma is immediate:

\begin{lemma}
Let $M' \subset M$ be any subgroup and $\Pi \subset \Sigma$. Then for any
$m \in H^D_\Pi$ the set $m + M'$ is contained in $H^D_\Pi$ if and only if $M'
\subset \bigcap_{\rho \in \Pi} \rho^\bot_M$.
\end{lemma}

\begin{proposition}\label{acycprop}
Let $B \subset \ksm$ be a monomial ideal, $D \in \Z^{\sigma(1)}$ and $\Pi \subseteq
\Sigma$, then:
\begin{enumerate}[(i)]
\item\label{acycpropi}
$\Gamma_B k[H^D_\Pi] =
\begin{cases}
k[H^D_\Pi] & \text{ if } \Pi \not \subset \Xi_B\\
0 & \text{ else}.
\end{cases}$
\item\label{acycpropii} The module $k[H^D_\Pi]$ is $\Gamma_B$-acyclic.
\end{enumerate}
\end{proposition}

\begin{proof}
Let first $\Pi \subset \Xi_B$, such that $\tau_\Pi \notin \supp(B)$. So there exists
an integral element $m$ in the relative interior of $\tau_\Pi^\bot \cap \sigma$ such
that the monomial $\chi(m)$ is contained in $B$. As the group $M' := \bigcap_{\rho
\in \Pi} \sigma^\bot_M$ contains $\tau_\Pi^\bot \cap \sigma_M$ (note that
$\tau_\Pi^\bot \cap \sigma_M = M' \cap \sigma_M$) and thus $m \in M'$, we have
$m + M' = M'$. This implies that any power of $\chi(m)$ is a nonzero divisor of the
module $k[H^D_\Pi]$, and moreover, multiplication by $\chi(m)$ even represents an
automorphism of the module $k[H^D_\Pi]$. So, any power of $\chi(m)$ also acts as an
automorphism on the local cohomology modules $H^i_B k[H^D_\Pi]$ for every $i \geq 0$.
But because the (ring theoretic) support of $k[H^D_\Pi]$ is contained in the support
of $B$, for every element $x \in H^i_B k[H^D_\Pi]$ there exists some $n > 0$ such
that $\chi(m)^n x = 0$, hence $H^i_B k[H^D_\Pi] = 0$ for every $i \geq 0$. So
for $\Pi \subset \Xi_B$, (\ref{acycpropi}) and (\ref{acycpropii}) are true.

Now consider the case $\Pi \not \subset \Xi_B$, i.e $\tau_\Pi \in \supp(B)$. In that
case, $B$ contains no monomial whose degree is contained in $\tau_\Pi^\bot \cap
\sigma_M$, and thus not in $M'$, so for every $x \in k[H^D_\Pi]$, there exists some
$n > 0$ such that $B^n x = 0$. So, the support of $k[H^D_\Pi]$ is contained in the
support of $B$, and thus $\Gamma_B k[H^D_\Pi] = k[H^D_\Pi]$ and $H^i_B k[H^D_\Pi]
= 0$ for every $i > 0$.
\end{proof}

\section{A $\Gamma_B$-Acyclic Resolution}\label{resolution}

Let $m \in M$ and define $\Sigma_m := \{\rho \in \Sigma \mid m \in H^D_\rho\}$.
Then we can consider the exact sequence $\mathbf{\Phi}_{\Sigma_m}$, respectively
$\mathbf{\Phi}_{\Sigma_m} \otimes_\Z k$. For every $-1 \leq i \leq \vert \Sigma_m
\vert$, we can identify the vector space at the $i$-th position of the complex
$\mathbf{\Phi}_{\Sigma_m} \otimes_\Z k$ as follows:
\begin{equation*}
\bigoplus_{\substack{\Pi \subset \Sigma_m\\ \vert \Pi \vert = i}} k \cdot \Pi \cong
\Big(\bigoplus_{\substack{\Pi \subset \Sigma_m\\ \vert \Pi \vert = i}} k[H^D_\Pi]
\Big)_m
\end{equation*}
With this identification for every $m \in M$, we obtain an exact sequence of vector
spaces:
\begin{equation*}
\mathbf{D}^\cdot: 
0 \rightarrow R^D \rightarrow k[M] \overset{\Phi_0}{\rightarrow}
\bigoplus_{\rho \in \Sigma} k[H^D_\rho] \overset{\Phi_1}{\rightarrow}
\bigoplus_{\substack{\Pi \subset \Sigma\\\vert \Sigma \vert = 2}} k[H^D_\Pi]
\overset{\Phi_2}{\rightarrow} \cdots \overset{\Phi_{n - 1}}{\rightarrow}
k[H^D_{\Sigma}] \rightarrow 0,
\end{equation*}
where the maps $\Phi_i$ are given by the direct sum of the cochain maps of the
complexes $\mathbf{\Phi}_{\Sigma_m}$ for every $m \in M$.

\begin{proposition}\label{acycres}
For any monomial ideal $B \subset \ksm$, the complex $\mathbf{D}^\cdot$ is a
$\Gamma_B$-acyclic resolution of $R^D$.
\end{proposition}

\begin{proof}
The exactness of $\mathbf{D}^\cdot$ already follows from its exactness as complex
of $k$-vector spaces. The $\Gamma_B$-acyclicity of the $k[H^D_\Pi]$ has already been
considered in the previous section, so it remains only to show that
$\mathbf{D}^\cdot$ indeed is a complex of $\ksm$-modules. For this, by $k$-linearity
it suffices to show that $\Phi_i \circ \chi(m) = \chi(m) \circ \Phi_i$ for any $-1
\leq i \leq \vert \Sigma \vert$ and $m \in \sigma_M$, where $\chi(m)$ is considered
as $k$-linear homomorphism. Consider any $m' \in M$, then we have inclusions
$\Sigma_{m' + m} \subset \Sigma_{m'} \subset \Sigma$, and we observe that
multiplication with $\chi(m)$ just results in the chain map
$\mathbf{\Phi}_{\Sigma_{m'}} \otimes_\Z k \twoheadrightarrow
\mathbf{\Phi}_{\Sigma_{m + m'}} \otimes_\Z k$, which yields the desired result.
\end{proof}

We define $\mathbf{D}^\cdot_{\Xi_B}$ as the subcomplex of $\mathbf{D}^\cdot$ which is
build of the terms
\begin{equation*}
\mathbf{D}_{\Xi_B}^i := \bigoplus_{\substack{\Gamma \subset \Pi \cap \Xi_B\\
\vert\Pi\vert = i}} k[H^D_\Gamma].
\end{equation*}
It is straightforward to see that degree-wise, for every $m \in M$, this complex
coincides with the complex $\mathbf{\Phi}_{\Sigma_m, \Xi_B} \otimes_\Z k$.

\begin{proposition}
$H^i_B R^D = H^i \big(\mathbf{D}^\cdot / \mathbf{D}^\cdot_{\Xi_B}\big)$ for every
$i \geq 0$.
\end{proposition}

\begin{proof}
Because $\mathbf{D}^\cdot$ is acyclic, we have that $H^i_B R^D = H^i\big(\Gamma_B(
\mathbf{D}^\cdot)\big)$. By proposition \ref{acycprop}, (\ref{acycpropii}),
$\Gamma_B(k[H^D_\Pi]) = 0$ if $\Pi \subset \Xi_B$ and $\Gamma_B k[H^D_\Pi] =
k[H^D_\Pi]$ if $\Pi \not \subset \Xi_B$, so the claim follows.
\end{proof}

\begin{corollary}\label{cohomcorollary}
For every $m \in M$ and every $i \geq 0$:
\begin{equation*}
\big(H^i_B R^D\big)_m = \tilde{H}^{i - 2}(\Xi_B \cap \Sigma_m; k).
\end{equation*}
\end{corollary}

\begin{proof}
Degreewise, for every $m \in M$, the complex $\mathbf{D}^\cdot /
\mathbf{D}^\cdot_{\Xi_B}$ is the augmented complex of relative cohomology with
coefficients in $k$ of the pair $(\Sigma_m, \Xi_B \cap \Sigma_m)$, shifted by 1.
So, in every degree,
we have an identification of the cohomology groups $\big(H^i R^D\big)_m =
\tilde{H}^{i - 1} (\Sigma_m, \Xi_B \cap \Sigma_m; k)$. Evaluating the long exact
cohomology sequence, using that $\tilde{H}^i(\Sigma_m; k) = 0$ for all $i \in \Z$, we
obtain that $\tilde{H}^{i - 1}(\Sigma_m, \Xi_B \cap \Sigma_m; k) \cong
\tilde{H}^{i - 2}(\Xi_B \cap \Sigma_m; k)$ for every $i \in \Z$.
\end{proof}

\section{The Chambers in $M_\R$ Determined by $R^D$}\label{chambersection}

In this section we assume the divisor $D = \sum_{\rho \in \Sigma} n_\rho D_\rho$ to
be fixed, except where explicitly stated otherwise.
For understanding the local cohomology modules $H_B^i R^D$, it is
important to know whether for some $\Pi \subset \Sigma$ a system of inequalities
\begin{align*}
\langle m, n(\rho) \rangle & < -n_\rho, \text{ for } \rho \in \Pi\\
\langle m, n(\rho) \rangle & \geq -n_\rho, \text{ for } \rho \in \Sigma \setminus \Pi
\end{align*}
has integral solutions or not.
For every $\rho \in \Sigma$, the linear equation $\langle m, n(\rho) \rangle =
-n(\rho)$ defines a hyperplane $H_\rho$ in $M$, and the set of hyperplanes $H_\rho$,
$\rho \in \Sigma$ forms a {\em hyperplane arrangement} in $M_\R$. The set of
inequalities above determines a {\em chamber} $C^{ss}_\Pi$ in the complement of this
hyperplane arrangement, i.e., if nonempty, the closure of the set of points $m$
fulfilling these inequalities form a polyhedron bounded by the hyperplanes $H_\rho$.
To be more precise, we define $C^s_\Pi$ to be the set of points fulfilling the
{\em strict} inequalities
\begin{align*}
\langle m, n(\rho) \rangle & < -n_\rho, \text{ for } \rho \in \Pi\\
\langle m, n(\rho) \rangle & > -n_\rho, \text{ for } \rho \in \Sigma \setminus \Pi
\end{align*}
and $C_\Pi$ analogously, but allowing equality in both types of equations. Note that
'ss' above stands for {\em semi-strict} inequalities. Moreover, note that we have
for simpler notation omitted any reference to the divisor $D$. The complement $M_\R
\setminus \bigcup_{\rho \in \Sigma} H_\rho$ then equals the set $\bigcup_{\Pi \subset
\Sigma} C^{s}_\Pi$. We have the following:

\begin{lemma}
The chambers $C^{ss}_\Pi$ and $C^s_\Pi$ are in one-to-one correspondence.
\end{lemma}

\begin{proof}
We show that if some point $m$ is contained in $C^s_\Pi$, then there exists a point
$m'$ which is contained in $C^{ss}_\Pi$, and vice versa. Let first $m \in C^s_\Pi$ for
some $\Pi \subset \Sigma$, then it is clearly contained in $C^{ss}_\Pi$. Now let $m
\in C^{ss}_\Pi$ for some $\Pi \subset \Sigma$. Let $\Gamma \subset \Sigma \setminus
\Pi$ given by precisely those $\rho$ such that $\langle m, n(\rho) \rangle = 0$. Now,
as the strict inequalities form an open condition, we can choose an
$\epsilon$-neighbourhood $U_\epsilon(m)$ in $M_\R$ such that for all points $x \in
U_\epsilon(m)$ the same strict inequalities hold as for $m$. Now, the inequalities
$\langle m, n(\rho) \rangle \geq 0$ for $\rho \in \Gamma$ determine a convex,
unbounded, polyhedron in
$M_\R$ which is not contained in a proper subspace of $M_\R$, and $m$ is located in
its boundary. Thus $U_\epsilon(m)$ must intersect the interior of this polyhedron,
and we can choose some $m'$ from the intersection of
$U_\epsilon$ and the interior of the polyhedron. Then it follows that $m' \in
C^s_\Pi$.
\end{proof}

To better understand the chambers type $C^{ss}_\Pi$, we have also to consider the
chambers of type $C_\Pi$, which are closed polyhedra in $M_\R$. Denote
$\sigma_\Pi$ the convex polyhedral cone generated over $\R_{\geq 0}$ by the lattice
vectors $-n(\rho)$ for $\rho \in \Pi$ and by $n(\rho)$ for $\rho \in \Sigma \setminus
\Pi$. Then every $C_\Pi$ can be written as Minkowski sum $P_\Pi + \check{\sigma}_\Pi$,
where $P_\Pi$ is a compact polyhedron and $\check{\sigma}_\Pi$ is the {\em dual} cone
of $\sigma_\Pi$. Our first observation
is that the polyhedra $C_\Pi$ do not have lineality spaces (see also \cite{Ziegler},
\S 1.5):

\begin{lemma}
The cones $\sigma_\Pi$ have dimension $d$ in $N_\R$.
\end{lemma}

\begin{proof}
Assume the vectors $-n(\rho)$ for $\rho \in \Pi$ and $n(\rho)$ for $\rho \in \Sigma
\setminus \Pi$ span a proper subvector space of $N_\R$, then also the vectors
$n(\rho)$, for $\rho \in \Sigma$, span a proper subvector space, but this is not
possible, since $\dim \sigma = d$.
\end{proof}

So, the cone $\check{\sigma}_\Pi$ can be identified with the {\em recession cone}
(sometimes also called the {\em characteristic cone}) of the polyhedron $C_\Pi$.
In general, $\check{\sigma}_\Pi$ is not $d$-dimensional. $\check{\sigma}_\Pi$
being $d$-dimensional is equivalent to that $\sigma_\Pi$ does not contain a nonzero
subvector space of $N_\R$. A general criterion for $\sigma_\Pi$ not containing a
nonzero subvector space is obtained by checking the intersection of the two cones
$\sigma^1_\Pi$ and $\sigma^2_\Pi$, where $\sigma_1$ is spanned over $\R_{\geq 0}$ by
$n(\rho)$, $\rho \in \Pi$ and $\sigma_2$ is spanned over $\R_{\geq 0}$ by $n(\rho)$,
$\rho \in \Sigma \setminus \Pi$.

\begin{lemma}
$\sigma_\Pi$ contains a nonzero subvector space if and only if $\sigma^1_\Pi \cap
\sigma^2_\Pi \neq \{0\}$.
\end{lemma}

\begin{proof}
Assume first that $\sigma_\Pi$ contains a nonzero subvector space $V$ and let $n \in
V$. Then we can write $n = \sum_{\rho \in \Pi} \alpha_\rho \big(-n(\rho)\big) +
\sum_{\rho \in \Sigma \setminus \Pi} \beta_\rho n(\rho)$, where $\alpha_\rho,
\beta_\rho \geq 0$. Because also $-n \in V$, we have $-n = \sum_{\rho \in \Pi}
\gamma_\rho \big(-n(\rho)\big) + \sum_{\rho \in \Sigma \setminus \Pi} \delta_\rho
n(\rho)$, for $\gamma_\rho, \delta_\rho \geq 0$. Summing up, we obtain
$n - n = 0 = \sum_{\rho \in \Pi} (\alpha_\rho + \gamma_\rho) \big(-n(\rho)\big) +
\sum_{\rho \in \Sigma \setminus \Pi} (\beta_\rho + \delta_\rho)n(\rho)$ where not
all of the $\alpha_\rho$, $\gamma_\rho$ and not all of the $\beta_\rho$, $\delta_\rho$
are zero, and thus the nonzero element $\sum_{\rho \in \Pi} (\alpha_\rho +
\gamma_\rho)n(\rho) = \sum_{\rho \in \Sigma \setminus \Pi} (\beta_\rho + \delta_\rho)
n(\rho)$ is contained in $\sigma^1_\Pi$ and $\sigma^2_\Pi$.
In the other direction, let $0 \neq n \in \sigma^1_\Pi \cap \sigma^2_\Pi$, i.e.
$n = \sum_{\rho \in \Pi} \alpha_\rho n(\rho) = \sum_{\rho \in \Sigma \setminus \Pi}
\beta_\rho n(\rho)$, where $\alpha_\rho, \beta_\rho \geq 0$ not all zero. Then $n \in
\sigma_\Pi$ and $-n = \sum_{\rho \in \Pi} \alpha_\rho \big(-n(\rho)\big) \in
\sigma_\Pi$ and thus $\sigma_\Pi$ contains the subvector space spanned by $n$.
\end{proof}

It would be nice if one could determine the dimension of the recession cone solely by
the combinatorics involved with the $\Pi$s, but in general this seems not to be
possible, as example \ref{noncombinatoriallyexample} below will show. However, at
least in the two extremal cases
there are some tools available. For the case of vanishing recession cones, the
corresponding chambers $C^s_\Pi$ are bounded, and these bounded chambers can at least
be counted by means of the {\em matroid of flats} associated to the hyperplane
arrangement $\bigcup_{\rho \in \Sigma} H_\rho$. For this, we refer to the book
\cite{BLSWZ}, chapter 4, in particular Corollary 4.6.8. The theory developed there
also yields a formula for counting all chambers $C^s_\Pi$.

The other extremal case is that of the recession cone $\check{\sigma}_\Pi$ having
dimension $d$. For these we can make use of the connection to local cohomology as
developed in the previous sections.

\begin{proposition}\label{dimpliescontractible}
Let $\check{\sigma}_\Pi$ be $d$-dimensional. Then either $\Pi = \Sigma$ or $\Pi \cap
\Xi$ is contractible.
\end{proposition}

\begin{proof}
First note that $\Pi = \Sigma$, then $\sigma_\Pi$ is just the negative cone of
$\sigma$, so we may assume that $\Pi \neq \Sigma$.
Let $\check{\sigma}_\Pi$ be $d$-dimensional and assume that $D = 0$, i.e. $R^D =
\ksm$. Then the corresponding hyperplane arrangement $\bigcup_{\rho \in \Sigma}
H_\rho$ is a central arrangement and the interior of $\check{\sigma}_\Pi$ coincides
with the chamber $C^s_\Pi$ in the complement of this arrangement, which therefore is
nonempty. Moreover, $C^s_\Pi$ has nonempty intersection with $M$. Choosing some $m
\in M \cap C^s_\Pi$, we can compute the local cohomology $H^i_{\mathfrak{m}} \ksm$ in
degree $m$, where $\mathfrak{m}$ is the maximal homogeneous ideal of $\ksm$. A
well-known result states that the ring $\ksm$ is
Cohen-Macaulay, and thus all these local cohomology modules for $0 \leq i < d$ vanish.
This in particular implies by Corollary \ref{cohomcorollary} that the reduced
cohomology groups $\tilde{H}^{i - 2}(\Pi \cap \Xi, k)$ vanish for $0 \leq i < d$.
So, because $\Pi \neq \Sigma$, this implies that $\Pi \cap \Xi$ is a contractible
topological space.
\end{proof}

\paragraph{The case $d = 3$.}

We show that in the case $d = 3$, there are possible only two types of recession cones
$\sigma_\Pi$, namely either $\sigma_\Pi$ is strictly convex, or $\sigma_\Pi
= N_\R$. For $d = 3$, the cell complex $\Xi$ is a topological 1-dimensional sphere
which can explicitly be realized as follows.
Choose a hyperplane $H$ of $N_\R$ such that $\sigma \cap H =: P$ is a bounded
polyhedron. Then the vertices of $P$ are given by $H \cap \rho$ and the facets of
$P$ coincide with the facets of $\sigma$ intersected with $H$. The set $\Xi$ then is
geometrically realized as the union of all $1$- and $2$-dimensional faces of $\sigma$
intersected with $H$, i.e., $\Xi$ has an explicit realization as the boundary of a
convex polytope in the plane
$H$. Moreover, $\Pi \neq \Sigma$ can be identified with a union of closed intervals
in $\Xi$. It is straightforward to see that the two cones $\sigma^1_\Pi$ and
$\sigma^2_\Pi$ intersect nontrivially if and only if the sets $P_1$, $P_2$ intersect,
where $P_1$ and $P_2$ are convex hulls of the points $\rho \cap H$, $\rho \in \Pi$ and
of the points $\rho \cap H$, $\rho \in \Sigma \setminus \Pi$, respectively.
We have the following

\begin{proposition}\label{threecontractible}
The cones $\sigma_\Pi^1$ and $\sigma_\Pi^2$ intersect nontrivially if and only if
$\Pi$ consists of more than one interval. In that case $\sigma_\Pi = N_\R$.
\end{proposition}

\begin{proof}
First note that, because $\Xi$ is a circle, $\Pi$ consists of as many intervals as
$\Sigma \setminus \Pi$. It follows from elementary geometric considerations that
in the case $\Pi$ consists of one interval, the polytopes $P_1$ and $P_2$ can not
intersect, and in the case where $\Pi$ consists of more than one interval, one can
choose vertices $p_1, p_2 \in P_1$ and $q_1, q_2 \in P_2$ such that the lines
$l_1 = p_1 + r(p_2 - p_1), r \in \R$, $l_2 = q_1 + s(q_2 - q_1), s \in \R$,
intersect in some point $a$ different from $p_1, p_2, q_1, q_2$. By arguments used
before, this implies that $\sigma_\Pi$ contains the subvector space spanned by $a$,
and because $a$ lies in the relative interior of the cone spanned by $q_1, q_2$,
$\sigma_\Pi$ also contains the 2-dimensional vector space spanned by $q_1, q_2$ over
$\R$. Moreover, the points $p_1$, $p_2$ are contained respectively on both, the
positive and the negative side of this vector space, so $\sigma_\Pi = N_\R$.
\end{proof}

Now we can prove:

\begin{theorem}
In the case $d = 3$, for $\Pi \subset \Sigma$, there are the following possibilities:
\begin{enumerate}[(i)]
\item $\tilde{H}^0(\Pi \cap \Xi, k) = 0$ and $C_\Pi$ has $d$-dimensional recession
cone $\sigma_\Pi$,
\item $\tilde{H}^0(\Pi \cap \Xi, k) \neq 0$ and $C^{ss}_\Pi$ is either bounded or
empty.
\end{enumerate}
\end{theorem}

\begin{proof}
We only observe that for $\Pi \neq \Sigma$, $\tilde{H}^0(\Pi \cap \xi, k) = 0$ is
equivalent to that $\Pi$ is contractible and apply propositions 
\ref{dimpliescontractible} and \ref{threecontractible}.
\end{proof}

We give the easiest example for the case $d = 3$:

\begin{example}\label{fourrayexample}
Let $\sigma$ be spanned over $\R_{\geq 0}$ by the primitive vectors $n_1 = (1, 0, 0)$,
$n_2 = (0, 1, 0)$, $n_3 = (-1, 1, 1)$, $n_4 = (0, 0, 1)$ and we consider the divisor
$D = -k D_2$ for some $k > 0$. For simplicity, we write here and in the
examples below $n_i$ instead of $n(\rho_i)$ and
$D_i$ instead of $D_{\rho_i}$. The hyperplane arrangement determined by $D$ realizes
precisely 15 nonempty chambers out of 16 possible choices $\Pi \subset \Sigma$, with
a unique bounded chamber for $\Pi = \{\rho_1, \rho_3\}$. We have $\dim \tilde{H}^0
(\Pi \cap \Xi, k) = 1$, and so the local cohomology module $H^2_{\mathfrak{p}_\sigma}
R^D$ does not vanish if $C^{ss}_\Pi \cap M \neq \emptyset$. This is the case for every
$k > 1$. For $k = 1$, the $C^{ss}_\Pi \cap M = \emptyset$, and thus $R^D$ is a
Cohen-Macaulay module. We obtain another Cohen-Macaulay module by analogous
considerations for $D = -D_1$, where this time the unique bound chamber is realized
for $\Pi = \{\rho_2, \rho_3\}$. Altogether, in this example there are three
isomorphism classes of maximal Cohen-Macaulay modules, represented by $\ksm$ itself,
$R^{-D_1}$, and $R^{-D_2}$.
\end{example}

We conclude that in the case $d = 3$, for computing local cohomology of the modules
$R^D$ it suffices to check the disconnected subsets $\Pi \subset \Sigma$. In general,
this may not be so simple, as we will see in the following examples. The first example
shows that even for topologically nontrivial $\Pi$, in general, the chamber
$C^{ss}_\Pi$ must not be bounded.

\begin{example}
Let $\sigma$ be spanned over $\R_{\geq 0}$ by the primitive vectors $n_1 =
(1, 0, 0, 0)$, $n_2 = (0, 1, 0, 0)$, $n_3 = (-1, 1, 1, 0)$, $n_4 = (0, 0, 1, 0)$,
$n_5 = (0, 0, 0, 1)$. This is the cone from the previous example extended by one ray
in four-dimensional direction. We choose $D = -k D_2$ for some $k > 0$, and it is
straightforward to see that, because $n_5$ is orthogonal to the other $n_i$, the
number of chambers is double the number of chambers of the previous example. Moreover,
it is easy to see that this time there are {\em no} bounded chambers, but still, for
$\Pi := \{\rho_1, \rho_4\}$, we have that $\Pi \cap \Xi$ consists of two points, and
hence $\dim \tilde{H}^0(\Pi \cap \Xi, k) = 1$
\end{example}

The next example shows that the contractibility of $\Pi$ does not imply that the
recession cone $\sigma_\Pi$ is strictly convex. Moreover, the example shows that the
strict convexity of $\sigma_\Pi$ can not depend on the combinatorics of $\Pi$ in a
simple way, but also depends on the concrete embedding of the cone $\sigma$ in $N_\R$.

\begin{example}\label{noncombinatoriallyexample}
Consider the four-dimensional cone $\sigma$ spanned by $n_1 = (0, 0, 0, 1)$, $n_2 =
(1, 0, 0, 1)$, $n_3 = (0, 1, 0, 1)$, $n_4 = (0, 0, 1, 1)$, $n_5 = (1, 1, 0, 1)$,
$n_6 = (1, 0, 1, 1)$, $n_7 = (0, 1, 1, 1)$, $n_8 = (1, 1, 1, 1)$, i.e. $\sigma$
is spanned over the three-dimensional unit cube shifted to the hyperplane $x_4 = 1$.
Set $\Pi := \{\rho_1, \rho_3, \rho_4, \rho_6\}$. Then $\Pi \cap \Xi$ is contractible
and we have $\tilde{H}^i(\Pi \cap \Xi, k) = 0$ for all $i$. But, we have $n_3 + n_6 =
n_2 + n_7$, so the cones $\sigma^1_\Pi$ and $\sigma^2_\Pi$ intersect, and the
recession cone $\check{\sigma}_\Pi$ is of dimension smaller than $4$.

Now consider the cone $\sigma'$ which is spanned by the same $n_1, \dots, n_8$ as
$\sigma$, except that $n_4$ and $n_6$ are replaced by $n_4' = (0, -1, 1, 1)$ and
$n_6' = (1, -1, 1, 1)$. $\sigma'$ is combinatorially equivalent to $\sigma$, but by
straightforward computation one finds that $(\sigma')^1$ and $(\sigma')^2$ do not
intersect, and thus $\check{\sigma}_\Pi$ is a $d$-dimensional recession cone of
$C_\Pi$ and thus $C^{ss}_\Pi$ is nonempty for every module $R^D$ over $k[\sigma'_M]$.
\end{example}

\section{Maximal Cohen-Macaulay Modules of Rank One}
\label{mcmsection}

By the results of section \ref{resolution}, the problem of classifying maximal
Cohen-Macaulay modules (MCMs) of rank one now has essentially become a problem of
integer programming. To see this more clearly, let us reformulate the results for
this case. For $R^D$ being an MCM is equivalent to that all local cohomology modules
$H^i_{\mathfrak{m}} R^D$ vanish for $i < d$, where $\mathfrak{m}$ is the maximal
homogeneous ideal of $\ksm$. This in particular is equivalent
to the vanishing of the cohomology groups $\tilde{H}^{i - 2}(\Sigma_m \cap \Xi, k)$
for every $i < d$ and every $m \in M$. Now, as we have seen in example
\ref{fourrayexample}, not every $\Pi$ such that $C^{ss}_\Pi$ is nonempty in $M_\R$
equals $\Sigma_m$ for some $m \in M$, i.e. not every $C^{ss}_\Pi$ has nonempty
intersection with $M$ although it is realized in the complement of the arrangement
$\bigcup_{\rho \in \Sigma} H_\rho$. So let us state the MCM-condition for $R^D$ as
a theorem:

\begin{theorem}
$R^D$ is an MCM if and only if for every $\Pi \subset \Sigma$ one of the following
two conditions holds:
\begin{enumerate}[(i)]
\item $\tilde{H}^{i - 2}( \Pi \cap \Xi, k) = 0$ for all $i < d$,
\item $\tilde{H}^{i - 2}( \Pi \cap \Xi, k) \neq 0$ for some $i < d$ and the chamber
$C^{ss}_\Pi$ is either empty or has empty intersection with $M$.
\end{enumerate}
\end{theorem}

We also state another, equivalent formulation:

\begin{theorem}\label{integertheorem}
$R^D$ is an MCM if and only if for every $\Pi \subset \Sigma$ one of the following two
conditions holds:
\begin{enumerate}[(i)]
\item $\tilde{H}^{i - 2}( \Pi \cap \Xi, k) = 0$ for all $i < d$,
\item $\tilde{H}^{i - 2}( \Pi \cap \Xi, k) \neq 0$ for some $i < d$ and the system of
inequalities
\begin{align*}
\langle m, n(\rho) \rangle & < - n_\rho \text{ for } \rho \in \Pi \\
\langle m, n(\rho) \rangle & \geq - n_\rho \text{ for } \rho \in \Sigma \setminus \Pi
\end{align*}
has no integral solution.
\end{enumerate}
\end{theorem}

One can relate the classification problem for MCMs of rank one to the problem of
understanding hyperplane arrangements in $M_\R$ induced by the hyperplanes $H_\rho$,
which are shifts of hyperplanes $\rho^\bot$ corresponding to some cone $\sigma \in
N_\R$ (such that, in particular, the hyperplanes $H_\rho$ are rational). If one fixes
the combinatorial type of the hyperplane arrangement $\bigcup_{\rho \in \Sigma}
H_\rho$, say, its matroid of flats, then, in how many ways can this hyperplane
arrangement be realized by shifting hyperplanes $H_\rho$, while keeping the
combinatorial
type, such that the cells $C^{ss}_\Pi$ with some nonvanishing cohomology group do not
intersect $M$?

\section{Singularity Sets}\label{singsets}

In order to actually proof that some module $R^D$ is an MCM, one effectively has to
check the
inequalities of theorem \ref{integertheorem} for nearly all possible sets $\Pi \subset
\Sigma$ which in general is a quite expensive task. In practice, however, it might be
a better strategy to check that some given $R^D$ is {\em not} an MCM. In the rest
of this paper we will collect some general results which can be helpful for this
purpose.

We introduce the notion of {\em singularity sets}; for the general theory of
singularity sets and their relation to local cohomology we refer to the book
\cite{SiuTrautmann}. For a variety $X$ over some algebraically closed
field $k$ and some coherent sheaf $\sh{F}$, the singularity sets of $\sh{F}$ are
defined for integers $i \geq 0$ as
\begin{equation*}
S_i(\sh{F}) := \{x \in X \mid \depth_{\sh{O}_{X, x}} \sh{F}_x \leq i\},
\end{equation*}
i.e. the set of points $x$ in $X$ such that the depth of the stalk $\sh{F}_x$ does not
exceed $i$.
The sets $S_i(\sh{F})$ are closed subsets of $X$ and every coherent sheaf defines a
filtration of $X$ by closed subsets $\emptyset \subset S_0(\sh{F}) \subset \cdots
\subset S_i(\sh{F}) \subset \cdots = X$. This filtration of course becomes stationary
for $i \geq \dim X$ with $S_{\dim X}(\sh{F}) = X$. We are only interested in the
situation where $X = U_\sigma$ is
an affine toric variety and $\sh{F} = \sh{O}(D)$, i.e. the sheafification of the
module $R^D$ over $U_\sigma$. Because $D$ is $T$-invariant, the depth of $\sh{O}(D)_x$
remains constant over every orbit $\orb(\tau) \subset U_\sigma$.
For $\tau \prec \sigma$, he restriction $\Gamma\big(U_\sigma, \sh{O}(D)\big)
\rightarrow \Gamma\big(U_\tau, \sh{O}(D)\big)$ corresponds to the localization
$R^D \rightarrow R^D_{\chi(m_\tau)}$, where $m_\tau$ is a lattice element from the
relative interior of the cone $\tau^\bot \cap \check{\sigma}$; in particular,
$k[\tau_M] = \ksm_{\chi(m_\tau)}$. Denote $\tau_{M_\tau} := \tau_M \cap M_\tau$, then
the semigroup $\tau_M$ splits into a cartesian product $\tau_{M} = \tau^\bot_M \times
\tau_{M_\tau}$. Correspondingly, the affine toric variety $U_\tau$ splits into the
cartesian
product $T^\tau \times U'_\tau$, where
$U'_\tau = \spec k[\tau_{M_\tau}]$ is an affine toric
variety of dimension $d - \codim \tau$. The corresponding projection $p: U_\tau
\twoheadrightarrow U'_\tau$ is a flat morphism. The following is a well-known fact on
equivariant sheaves or $M$-graded modules, respectively, which we present without
proof.

\begin{proposition}
Every $T$-equivariant coherent sheaf $\sh{E}$ over $U_\tau$ is isomorphic to
$p^*\sh{E}'$ for some $T^\sigma$-equivariant sheaf $\sh{E}'$ over $U'_\tau$.
Equivalently, every finitely generated $M$-graded $k[\tau_M]$-module $E$ is isomorphic
to $E' \otimes_{k[\tau'_M]} k[\tau_M]$ for some $M_\tau$-graded $k[\tau'_M]$-module
$E'$.
\end{proposition}

In particular, $\sh{O}(D) \cong p^* \sh{O}(D')$, where $D' = \sum_{\rho \in \tau(1)}
n_\rho D_\rho$ is a $T^\sigma$-invariant divisor on $U'_\tau$ and $\sh{O}(D')$ is the
sheafification of the $M$-graded $k[\tau_{M_\tau}]$-module $R^{D'}$. We obtain:

\begin{lemma}
For every $x \in U_\tau$, we have $\depth_{\sh{O}_{U_\tau, x}} \sh{O}(D)_x =
\depth_{\sh{O}_{U'_\tau, p(x)}} \sh{O}(D')_{p(x)} + \codim \tau$.
\end{lemma}

\begin{proof}
As the morphism $p$ is flat and thus local, we can apply \cite{Matsumura}, Thm. 23.3
and obtain $\depth_{k[\tau_M]_x} R^D_x = \depth_{k[\tau_M]_x} R^{D'}_{p(x)}
\otimes_{k[\tau_{M_\tau}]_{p(x)}} k[\tau_M]_x = \depth_{k[\tau_{M_\tau}]_{p(x)}}
R^{D'}_{p(x)} + \depth_{k[\tau_M]_x} \big(k[\tau_M] / \mathfrak{p}_\tau \big)_x$ for
every point $x \in U_\tau$, where here $\mathfrak{p}_\tau$ denotes the maximal
homogeneous ideal of $k[\tau_M]$.
\end{proof}

With help of this lemma, we set:

\begin{definition}
Let $i \geq 0$, then we set
\begin{equation*}
S_i := \{\tau \prec \sigma \mid \depth_{\ksm_x} R^D_x \leq i - \codim \tau \text{ for
some point $x \in \orb(\tau)$}\}.
\end{equation*}
\end{definition}

In this definition, we have omitted any explicit reference to $D$ for clearer
notation. Note that $S_i$ is star-closed, and by the discussion above,
$S_i\big(\sh{O}(D)\big)$ is equal to $\bigcup_{\tau \in S_i} \orb(\tau)$ for all
$i \geq 0$. Now observe:

\begin{lemma}
$R^D$ is MCM if and only if $S_i = \emptyset$ for $0 \leq i < d$.
\end{lemma}

Denote $\Xi^\tau := \{\Pi \prec \tau(1) \mid \tau_\Pi \notin \stern(\tau)\}$, where we
consider $\Xi^\tau$ as a subcomplex of the simplex of $\tau$.
The following lemma is immediate:

\begin{lemma}\label{tinylemma}
Let $\tau \prec \sigma$ and $\Pi \subset \tau(1)$, then $\Pi \cap \Xi^\tau = \Pi \cap
\Xi_\tau$.
\end{lemma}

For any subset $\Pi \subset \tau(1)$, the splitting $M \cong \tau_M^\bot \times M_\tau$
is compatible with linear
inequalities $\langle m, n(\rho) \rangle < -n_\rho$, for $\rho \in \Pi$,
respectively $\langle m, n(\rho) \rangle \geq -n_\rho$ for $\rho \in \tau(1) \setminus
\Pi$, in the sense that some $m \in M$ fulfills these inequalities if and only if
every $m'$ with $m' - m \in \tau_M^\bot$ fulfills these inequalities.
The following is the main result of this section:

\begin{theorem}
$\tau \in S_k$ if and only if $\tilde{H}^{i - 2}(\Pi \cap \Xi_{\tau}, k) \neq 0$ for
some $i \leq k - \codim \tau$ and some subset $\Pi \subset \tau(1)$ and there
exists an integral solution to the system of inequalities
\begin{align*}
\langle m, n(\rho) \rangle & < - n_\rho \text{ for } \rho \in \Pi \\
\langle m, n(\rho) \rangle & \geq - n_\rho \text{ for } \rho \in \tau(1) \setminus \Pi
\end{align*}
in $M$ or equivalently, in $M_\tau$.
\end{theorem}

\begin{proof}
We have $\tau \in S_k$ if and only if $H^i_{\mathfrak{p}_\tau} R^{D'} \neq 0$ for some
$i \leq k - \codim \tau$, where $\mathfrak{p}_\tau$ is the maximal homogeneous ideal
of $k[\tau_{M_\tau}]$. This in turn is equivalent to that there exists some $m \in
M_\tau$ and some $i \leq k - \codim \tau$ such that $\tilde{H}^i(\Sigma^\tau_m \cap
\Xi^\tau, k) \neq 0$, where $\Sigma^\tau_m = \{\rho \in \tau(1) \mid \langle m, n(\rho)
\rangle < -n_\rho\}$. Moreover, by lemma \ref{tinylemma}, $\Sigma^\tau_m \cap \Xi^\tau
= \Sigma^\tau_m \cap \Xi_\tau$.
\end{proof}

The theorem can help to reduce the number of inequalities one has to check in order to
determine the sets $S_i$. However, in the case $d = 3$, this is not of much help.

\begin{proposition}[\cite{SiuTrautmann}, Corollary 1.21]
Let $X$ be an irreducible variety of dimension $d$ and let $\sh{F}$ be a coherent sheaf
on $X$, then $\sh{F}$ is reflexive if and only if $\dim S_i(\sh{F}) \leq i - 2$ for all
$i < d$.
\end{proposition}

For $d = 3$ this implies that $S_2$ is either $\{\sigma\}$ or empty.

\newcommand{\etalchar}[1]{$^{#1}$}

\end{document}